\documentclass{amsart}
\usepackage{graphicx}
\usepackage{colordvi}

\newtheorem{thm}{Theorem}[section]
\newtheorem{cor}[thm]{Corollary}
\newtheorem{lem}[thm]{Lemma}

\newtheorem{rem}[thm]{Remark}

\newcommand{\abs}[1]{\left\vert#1\right\vert}
\newcommand{\set}[1]{\left\{#1\right\}}
\newcommand{\Real}{\mathbb R}

\newcommand{\pfrac}[2]{\frac{\partial #1}{\partial #2}}

\title[Ricci flow on surfaces]{Normalized Ricci flow on nonparabolic surfaces
} 
\author{Hao Yin}

\begin{document}                                                                                   
\maketitle
\begin{abstract}
This paper studies normalized Ricci flow on a nonparabolic surface, whose scalar curvature is asymptotically $-1$ in an integral sense. By a method initiated by R. Hamilton, the flow is shown to converge to a metric of constant scalar curvature $-1$. A relative estimate of Green's function is proved as a tool.
\end{abstract}

\section{Introduction}\label{S_Intro}

Let $(M,g)$ be a Riemannian manifold of dimension 2. The normalized ricci flow is
\begin{equation*}
\pfrac{g_{ij}}{t}=(r-R)g_{ij},
\end{equation*}
where $R$ is the scalar curvature and $r$ is some constant. For compact surface, $r$ is the average of scalar curvature. In this case, Hamilton \cite{Ha} and Chow \cite{Chow} proved the normalized Ricci flow from any initial metric will exist for all time and converge to a metric of constant curvature. It's therefore nature to ask if such result holds for non-compact surfaces. Recently, a preprint of Ji and Sesum \cite{Se} generalized the above result to complete surfaces with logarithmic ends. Such surfaces have infinities like hyperbolic cusps. In particular, they have finite volume, therefore are parabolic, in the sense that there exists no positive Green's function. One of their result shows that the normalized Ricci flow from such a metric will exist for all time and converge to hyperbolic metric. In this paper, we study nonparabolic complete surfaces, i.e. surfaces admitting positive Green's function. In contrast to \cite{Se}, such surfaces have at least one nonparabolic end and have infinite volume. For a discussion of parabolic and nonparabolic ends and their geometric characterization, see Li's survey paper \cite{LiSurvey}.

Here we choose $r=-1$ because if the flow converges, the limit metric will be of constant curvature $r$. Since we are considering noncompact surfaces, $r$ can't be positive. If $r=0$, the limit will be flat $\Real^2$ or its quotient. However, it's well known that these flat surfaces are parabolic. On the other hand, whether a surface is parabolic or nonparabolic is invariant under quasi-isometries. Since if the normalized Ricci flow converges, then the limit metric will be quasi-isometric to the initial one, we know $r$ can't be zero.(For the definition of quasi-isometry, see also \cite{LiSurvey}.) If $r<0$, we can always assume $r=-1$ by a scaling.

The main result of this paper is
\begin{thm}
\label{thm_main}
Let $(M,g)$ be a  nonparabolic surface with bounded curvature. If the infinity is close to a hyperbolic metric in the sense that
\begin{equation*}
\int_M \abs{R+1} dV<+\infty.
\end{equation*}
Then, the normalized Ricci flow will converge to a  metric of constant scalar curvature $-1$.
\end{thm}

As in \cite{Se}, we try to apply the above result to prove results along the line of Uniformization theorem. That amounts to prove the existence of a complete hyperbolic metric within a given conformal class of a noncompact surface. In \cite{Se}, the authors proved that there is a uniformization theorem for Riemann surfaces obtained from compact Riemann surface by removing finitely many points and remarked that similar result should be true for Riemann surfaces obtained from compact ones by removing finitely many disjoint disks and points. Our theorem can be used to prove the same result in the case there is at least one disk removed. In fact, we will give a unified proof, which includes and simplifies the proof of \cite{Se}. Precisely, we will show

\begin{cor}\label{cor}
Let $M$ be a Riemann surface obtained from compact Riemann surface by removing finitely many disjoint disks and/or points. If no disk is removed, then we further assume that the Euler number of $M$ is less than zero. Then there exists on $M$ a complete hyperbolic metric compatible with the conformal structure.
\end{cor}

The proof of Theorem \ref{thm_main} is along the same line as \cite{Se}. The method was initiated by Hamilton in \cite{Ha}. There, Hamilton considered only compact case. for the purpose of generalizing this method to complete case, we need to overcome some analytic difficulties. Precisely, one need to solve Poisson equations and obtain estimates for the solutions, for all $t$. Those growth estimates for the solution are needed to apply the maximum principle. As for the maximum principle, there are many versions of maximum principle on complete manifolds. Since we will be working on complete manifold with a changing metric, the closest version for our need is in \cite{Tam}. We still need a little modification. 

\begin{thm}
	\label{thm_maximum}
	Suppose $g(t)$ is a smooth family of complete metrics defined on $M$, $0\leq t\leq T$ with Ricci curvature bounded from below and $\abs{\pfrac{}{t}g}\leq C$ on $M\times [0,T]$. Suppose $f(x,t)$ is a smooth function defined on $M\times [0,T]$ such that
	\begin{equation*}
		\triangle^t f-\pfrac{f}{t}\geq 0
	\end{equation*}
	whenever $f(x,t)>0$
	and
	\begin{equation}
		\label{eqn_growth}
		\int_0^T \int_M \exp(-a r_t^2(o,x))f_+^2(x,t)dV_t<\infty
	\end{equation}
	for some $a>0$. If $f(x,0)\leq 0$ for all $x\in M$, then $f\leq 0$ on $M\times [0,T]$.
\end{thm}
Although there is no detail in \cite{Tam}, one can prove it using the method of Ecker and Huisken in \cite{Huisken} and Ni and Tam in \cite{NT}.

To solve the Poisson equation $\triangle u=R+1$ for $t=0$. We use a result of Ni\cite{Ni}, See Theorem \ref{thm_Ni}. That's the reason why we assume $\int_M \abs{R+1}dV< +\infty$. Moreover, we prove a growth estimate of the solution under the further assumption that Ricci curvature bounded from blow. This result is true for all dimensions. For the growth estimate, an estimate of Green's function is proved under the assumption that Ricci curvature bounded from below. This estimate may be of independent interests, see the discussion in Section \ref{S_Green}.

Instead of solving $\triangle_t u(x,t)=R(x,t)+1$ for later $t$. We solve an evolution equation for $u$. Thanks to the recent preprint of Chau, Tam and Yu \cite{Tam}, we can solve this evolution equation with a changing metric. Following a method in \cite{NiTam}, we show that $u$, $\abs{\nabla u}$ and $\triangle u$ satisfy the growth estimate like in equation (\ref{eqn_growth}). With these preperation, we proceed to show that $u(x,t)$ is indeed the potential functions we need. Now the Theorem \ref{thm_main} follows from the approach of Hamilton and repeated use of Theorem \ref{thm_maximum}.

The paper is organized as follows: In Section \ref{S_Green}, we prove the crucial estimate of Green's function needed for the growth estimate. In Section \ref{S_Poisson}, we solve the Poisson equation and prove the relevant growth estimates. In the last section, we prove Theorem \ref{thm_main} and discuss results related to Uniformization theorem.

\section{An estimate of Green's function}\label{S_Green}

In this section we prove that
\begin{thm}
	\label{thm_Green}
	Let $(M,g)$ be a complete noncompact manifold with Ricci curvature bounded from below by $-K$. Assume that $M$ admits a positive Green's function $G(x,y)$. Let $x_0$ be a fixed point in $M$. Then there exists constant $A>0$ and $B>0$, which may depend on $M$ and $x_0$, so that
	\begin{equation*}
		\int_{\{G(x,y)>e^{Ar(y,x_0)}\}}G(x,y)dx\leq Be^{Ar(y,x_0)},
	\end{equation*}
	where $r(y,x_0)$ is the distance from $y$ to $x_0$.
\end{thm}

\begin{rem}
	It's impossible to get an estimate of this kind with constant depending only on $K$. Considering a family of nonparabolic manifolds $M_i$, which are becoming less and less 'nonparabolic', i.e. their infinities are closing up. For any $A,B>0$, there exists $M_i$ and some $x_i\in M_i$ such that
	\begin{equation*}
		\int_{\{G_i(x,x_i)>A\}}G_i(x,x_i)dx>B.
	\end{equation*}
	See \cite{LT}.
\end{rem}

\begin{rem}
	To the best of the author's knowledge, known estimates on Green's function in terms of volume of balls require Ricci curvature to be non-negative, See \cite{LY}. There could be one estimate of such type for Ricci curvature bounded from below, in light of \cite{Tam}. If so, our relative estimate should be a corollary. The following proof is a direct one.
\end{rem}

We begin with a lemma,
\begin{lem}
	\label{lem_Green}
	There is a constant $C$ depending only on $K$ and the dimension, such that if Ricci curvature on $B(x,1)$ is bounded from below by $-K$ and $G(x,y)$ is the Dirichlet Green's function on $B(x,1)$, then
	\begin{equation*}
		\int_{B(x,1)}G(x,y)dy<C.
	\end{equation*}
\end{lem}

\begin{proof}
	Let $H(x,y,t)$ be the Dirichlet heat kernel of $B(x,1)$. It's easy to see
	\begin{equation*}
		\int_{B(x,1)} H(x,y,t)dy\leq 1,
	\end{equation*}
	for all $t>0$.
	
	Now we prove that $H(x,y,2)$ is bounded from above. The proof is Moser iteration, which has appeared several times. Here we follow computations in \cite{Zhang}. Since we have Dirichlet boundary condition, we don't need cut off function of space.

	Let $0<\tau<2$ and $0<\delta\leq 1/2$ be some positive constants, $\sigma_k=(1-(1/2)^k \delta )\tau$ and $\eta_i$ be smooth function on $[0,\infty)$ such that 1) $\eta_i=0$ on $[0,\sigma_i]$, 2) $\eta_i=1$ on $[\sigma_{i+1},\infty)$ and 3) $\eta_i^\prime\leq 2^{i+3}(\delta\tau)^{-1}$. Let $p_i=(1+\frac{2}{n})^i$.
	Since $H$ is a solution to the heat equation, it's easy to know $H^p$ is a subsolution to the heat equation for $p>1$.
	\begin{equation*}
		(\pfrac{}{t}-\triangle_y) H^p(x,y,t)\leq 0.
	\end{equation*}
	Multiply by $\eta_i^2 H^{p_i}$ and integrate
	\begin{equation*}
		\int_{\sigma_i}^T\int_{B(x,1)} \eta_i^2H^{p_i}(\pfrac{}{t}-\triangle_y)H^{p_i}dydt\leq 0.
	\end{equation*}
	Routine computation gives
	\begin{equation*}\label{equ_kaka}
		\int_{\sigma_{i+1}}^T\int_{B(x,1)} \abs{\nabla_y H^{p_i}}^2 dy dt+ \frac{1}{2}\int_{B(x,1)}H^{2p_i}(x,y,T)dy \leq 2^{i+3}(\tau\delta)^{-1}\int_{\sigma_i}^T\int_{B(x,1)}H^{2p_i}dydt.
	\end{equation*}
	The sobolev inequality in \cite{Sobolev} implies
	\begin{equation*}
		\left(\int_{B(x,1)}(H^{p_i})^{\frac{2n}{n-2}}dy\right)^{\frac{n-2}{n}} \leq C V^{-2/n}  \int_{B(x,1)}\abs{\nabla_y H^{p_i}}^2 + H^{2p_i}dy,
	\end{equation*}
	where $V$ is the volume of $B(x,1)$.
	By H\"{o}lder inequality,
	\begin{eqnarray*}
		\int_{B(x,1)} H^{2p_{i+1}}dy &\leq& \left(\int_{B(x,1)}(H^{p_i})^{\frac{2n}{n-2}}dy\right)^{\frac{n-2}{n}} (\int_{B(x,1)} H^{2p_i}dy)^{2/n} \\
		&\leq& (C V^{-2/n}  \int_{B(x,1)}\abs{\nabla_y H^{p_i}}^2 + H^{2p_i}dy)(\int_{B(x,1)}H^{2p_i}dy)^{2/n}.
	\end{eqnarray*}
	By (\ref{equ_kaka}), integrate over time
	\begin{eqnarray*}
		\int_{\sigma_{i+1}}^2 \int_{B(x,1)} H^{2p_{i+1}}dydt &\leq& CV^{-2/n} c_0^{i+3} (\sigma\tau)^{-(1+2/n)} (\int_{\sigma_i} ^2 \int_{B(x,1)} H^{2p_i}dydt)^{1+\frac{2}{n}},
	\end{eqnarray*}
	where $c_0=2^{1+2/n}$.
	A standard Moser iteration gives
	\begin{equation*}\label{eqn_hehe}
		\max_{t\in [\tau,2]}\max_{y\in B(x,1)} H^2(x,y,t)\leq CV^{-1}(\sigma\tau)^{-\frac{n+2}{2}}\int_{(1-\delta)\tau}^2\int_{B(x,1)} H^2(x,y,t)dydt.
	\end{equation*}
	An iteration process as given in \cite{LS} implies the $L^1$ mean value inequality. In particular,
	\begin{equation*}
		\max_{y\in B(x,1)} H(x,y,2)\leq CV^{-1} \int_1^2\int_{B(x,1)}H(x,y,t)dy dt\leq C V^{-1}.
	\end{equation*}
	Hence,
	\begin{equation*}
		\int_{B(x,1)}H^2(x,y,2)dy\leq C V^{-1}.
	\end{equation*}
	Due to a Poincar\'{e} inequality in \cite{LS},
	\begin{eqnarray*}
		\frac{d}{dt}\int_{B(x,1)}H^2(x,y,t)dy &=& \int_{B(x,1)}2H\triangle_y H dy \\
		&=& -\int_{B(x,1)}\abs{\nabla_y H}^2 dy \\
		&\leq& -C \int_{B(x,1)}H^2(x,y,t)dy.
	\end{eqnarray*}
	This differential inequality implies
	\begin{eqnarray*}
		\int_{B(x,1)}H^2(x,y,t)dy &\leq & \int_{B(x,1)}H^2(x,y,2)dy \times e^{-C(t-2)} \\
		&\leq& CV^{-1}e^{-C(t-2)}.
	\end{eqnarray*}
	H\"{o}lder inequality shows
	\begin{equation*}
		\int_{B(x,1)}H(x,y,t)\leq V(B(x,1))\int_{B(x,1)}H^2(x,y,t)dy\leq Ce^{-C(t-2)},
	\end{equation*}
	for $t\geq 2$.
	The lemma follows from 
	\begin{equation*}
		\int_{B(x,1)}G(x,y)dy= \int_{0}^\infty\int_{B(x,1)}H(x,y,t)dydt.
	\end{equation*}
\end{proof}

Now let's turn to the proof of Theorem \ref{thm_Green}.

\begin{proof}
The key tool in the proof is Gradient estimate for harmonic function. Recall that if $u$ is a positive harmonic function on $B(x,2R)$, then
\begin{equation*}
\sup_{B(x,R)}\abs{\nabla \log u(x)}^2\leq C_1K+C_2R^{-2}
\end{equation*}

This is to say outside $B(x,0.1)$, the Green function as a function of $y$ decays or increases at most exponentially with a factor $\sqrt{C_1K+100C_2}$.

(1) Consider $G(x_0,y)$, Set
\begin{equation*}
p=\max_{y\in \partial B(x_0,1)} G(x_0,y).
\end{equation*}
As pointed out in Li and Tam, in the paper constructing Green function, $G(x_0,y)\leq p$ for $y\notin B(x_0,1)$. Since the Green function is symmetric, for any point $y$ far out in the infinity, $G(y,x_0)\leq p$.

(2) If the theorem is not true, then for any big $A$ and $B$, there is a point $y$ (far away) so that
\begin{equation*}
\int_{\{G(x,y)>e^{Ar(y,x_0)}\}}G(x,y)> Be^{Ar(y,x_0)}.
\end{equation*}
We will derive a contradiction with (1).

Claim: $\{x| G(x,y)>e^{Ar}\}\subset B(y,1)$ is not true.

If true, then consider the Dirichlet Green function $G_1(z,y)$ on $B(y,1)$. It's well known that $G(z,y)-G_1(z,y)$ is a harmonic function. Notice that this harmonic function has boundary value less than $e^{Ar}$. Therefore, its integration on $B(y,1)$ is less than $e^{Ar}\times Vol(B(y,1))$. Since we assume Ricci lower bound, $Vol(B(y,1))$ is less than a universal constant depending on $K$.

Therefore,
\begin{eqnarray*}
	\int_{\{G(x,y)>e^{Ar}\}} G(x,y)dx &\leq & \int_{B(y,1)} G(x,y)dx \\
	&\leq & Vol(B(y,1))\times e^{Ar} + \int_{B(y,1)} G_1(x,y) dx \\
	&\leq & C(K,n)\times e^{Ar},
\end{eqnarray*}
where we used Lemma \ref{lem_Green} for the last inequality. If we choose $B$ to be any number larger than $C(K,n)$ in the above equation, then the choice of $y$ gives an contradiction and implies that the claim is true.

(3) There is a $z\in \{x| G(x,y)>e^{Ar}\}$ so that $d(z,y)=1$ because the set $\{G(x,y)>e^{Ar}\}$ is connected. This follows from the maximum principle and the construction of Green's function.

(3.1) If $\abs{d(y,x_0)-d(z,x_0)}<0.3$, then

Let $\sigma$ be the minimal geodesic connecting $z$ and $x_0$.

Claim: the nearest distance from $y$ to $\sigma$ is no less than 0.1.

If not, let $w$ be the point in $\sigma$ such that $d(y,w)<0.1$. Since $d(y,z)>1$, we know
\begin{equation*}
d(w,z)>0.9
\end{equation*}
Now, $w$ is on the minimal geodesic from $z$ to $x_0$, so
\begin{equation*}
d(w,x_0)\leq d(z,x_0)-0.9
\end{equation*}
then
\begin{equation*}
d(y,x_0)<d(w,x_0)+d(y,w)< d(z,x_0)-0.8
\end{equation*}
This is a contradiction , so the claim is true.

We can use the gradient estimate along the segment $\sigma$. (Notice that $d(z,x_0)<r(x,x_0)+1$)
\begin{equation*}
G(y,x_0)>\frac{G(y,z)}{exp(\sqrt{C_1K+100C_2}(r+1))}
\end{equation*}
This is a contradiction if we choose $A>> \sqrt{C_1K+100C_2}$.

(3.2) If $d(z,x_0)\leq d(y,x_0)-0.3$, then

The distance from $y$ to the minimal geodesic connecting $z$ and $x_0$ will be larger than $0.1$. The above argument gives a contradiction.

(3.3) If $d(z,x_0)\geq d(y,x_0)+0.3$, then

Since $G(z,y)>e^{Ar}$, we move the center to $z$, by symmetry of Green function. $G(y,z)>e^{Ar(y,x_0))}>e^{A^\prime r(z,x_0)}$. This is case (3.2). We get a contradiction at $z$.

This finishes the proof of estimate of Green function.
\end{proof}

\section{Poisson equations $\triangle u=R+1$}\label{S_Poisson}

This section is divided into two parts. The first part solves the Poisson equation for $t=0$. The second part solves for $t>0$ before the maximum time using an indirect way. 

First,  we use Theorem \ref{thm_Green} to obtain an growth estimate of the solution of the Poisson equation $\triangle u=R+1$ for $t=0$. The existence part without curvature restriction and boundedness of $f$ of the following theorem is due to Lei Ni in \cite{Ni}. 

\begin{thm}
\label{thm_Ni}
Let $M$ be a complete nonparabolic manifold with Ricci curvature bounded from below by $-K$. For non-negative bounded continuous function $f$ the Poisson equation
\begin{equation*}
	\triangle u=-f
\end{equation*}
has a non-negative solution $u\in W^{2,n}_{loc}(M)\cap C^{1,\alpha}_{loc}(M)(0<\alpha<1)$ if $f\in L^1(M)$. Moreover, for any fixed $x_0\in M$, there exists $A>0$ and $C>0$ such that
\begin{equation*}
	u(x)\leq C e^{Ar(x,x_0)}.
\end{equation*}
\end{thm}

\begin{proof}
	Let $G(x,y)$ be the positive Green's function. 
	\begin{eqnarray*}
		\int_M G(x,y)f(y)dy &=& \int_{\{G(x,y)\leq e^{Ar(x,x_0)}\}} G(x,y)f(y) dy \\
&&+\int_{\{G(x,y)> e^{Ar(x,x_0)}\}} G(x,y)f(y)dy \\
		&\leq& Ce^{Ar(x,x_0)} .
	\end{eqnarray*}
	For the first term, we use the assumption that $f$ is integrable, for the second term, we use the boundedness of $f$ and the Theorem \ref{thm_Green}. The estimate above shows the Poisson equation is solvable with the required estimate.
\end{proof}

\begin{cor}
Let $M$ be a surface satisfying the assumptions in Theorem \ref{thm_main}. There exists a solution $u_0$ to the equation $\triangle u_0=R(x)+1$ satisfying
\begin{equation*}
	\int_M \exp(-a r^2(x,x_0))u_0^2(x) dV<\infty
\end{equation*}
and
\begin{equation*}
	\int_M \exp (-b r^2(x,x_0)) \abs{\nabla u_0}^2(x) dV<\infty
\end{equation*}
where $a$ and $b$ are some positive constants.
\end{cor}
\begin{proof}
Solve the Poisson equation for the positive part and the negative part of $R+1$ respectively. Then subtract the solutions. The first integral estimate follows from the pointwise growth estimate and volume comparison. 

Let $R>1$. Choose a cut-off function $\varphi$ such that
\begin{equation*}
	\varphi(x)=\left\{
	\begin{array}{ll}
		1 &  x\in B(x_0,R) \\
		0 & x\notin B(x_0,2R)
	\end{array}
	\right.
\end{equation*}
and
\begin{equation*}
	\abs{\nabla \varphi}^2\leq C_1\varphi.
\end{equation*}
Multiply the equation by $\varphi u_0$ and integrate over $M$,
\begin{equation*}
	\int_M \varphi u_0 \triangle u_0 dV =\int_{M} (R+1) \varphi u_0 dV,
\end{equation*}
which implies
\begin{equation*}
	\int_M \varphi\abs{\nabla u_0}^2 dV +\int_M u_0 \nabla \varphi\cdot \nabla u_0 dV =-\int_M (R+1) \varphi u_0dV.
\end{equation*}
Hence
\begin{eqnarray*}
	\int_M (\varphi-\frac{\abs{\nabla \varphi}^2}{2C_1})\abs{\nabla u_0}^2 dV &\leq& C\int_{B(x_0,2R)}u_0^2dV + C\int_{B(x_0,2R)} \abs{u_0} dV. \\
	&\leq& C\int_{B(x_0,2R)}u_0^2 dV +C Vol(B(x_0,2R)).
\end{eqnarray*}
From the integration estimate of $u_0$, 
\begin{equation*}
	\int_{B(x_0,2R)}u_0^2 dV\leq Ce^{4aR^2}.
\end{equation*}
By choice of $\varphi$,
\begin{equation*}
	\int_{B(x_0,R)} \abs{\nabla u_0}^2 dV\leq Ce^{\tilde{a}R^2}.
\end{equation*}
From here, it's not difficult to see the estimate we need.
\end{proof}

Now let's look at the case of $t>0$. In fact, it's not difficult to show the above method can be used for $t>0$. This amounts to show that $M$ is still nonparabolic for $t>0$ and $\int_M \abs{R+1}dV$ is still finite. The first claim is trivial and the second follows from the evolution equation and maximum principle. Assume the solutions are $u(t)$. We have trouble in deriving the evolution equation for $u(t)$, due to the possible existence of nontrivial harmonic functions. This explains why we use the following indirect way.

\begin{lem}
	Assume the normalized Ricci flow exists for $t\in [0,T_{max})$. The following equation has a solution $u(x,t) \quad (0\leq t<T_{max})$ with initial value $u_0$, 
\begin{equation*}
	\pfrac{u}{t}=\triangle u -u,
\end{equation*}
where $\triangle$ is the Laplace operator of metric $g(t)$. Moreover, there exists $a>0$ depending on $T$ such that for any $T<T_{max}$
\begin{equation*}
\int_0^T\int_M \exp(-a r^2(x,x_0)) u^2(x,t)dV_t <\infty.
\end{equation*}
Similar estimates hold for $\abs{\nabla u}$ and $\triangle u$ with different constants.
\end{lem}

\begin{rem}
	Since $g(0)$ and $g(t)$ are equivalent up to a constant depending on $T$, it doesn't matter whether we estimate $\nabla u$ or $\nabla_t u$ and whether we use $r$ to stand for distance at $g(0)$ or $g(t)$ if $t\in [0,T]$.
\end{rem}
\begin{proof}
	In \cite{Tam}, the authors considered a class of evolution equation with changing metric. $\pfrac{u}{t}=\triangle u -u$ with the underling metric evolving by normalized Ricci flow is in this class. They proved, among other things, that the fundamental solution $\mathcal{Z}(x,t;y,s)$ has a Gaussian upper bound, i.e\begin{equation*}
\mathcal{Z}(x,t;y,x)\leq \frac{C}{V_x(\sqrt{t-s})}e^{-\frac{r^2(x,y)}{D(t-s)}}.
\end{equation*}
These constants depends on the solution of normalized Ricci flow and $T$. See Corollary 5.2 in \cite{Tam}. For simplicity, denote $\mathcal{Z}(x,t;y,0)$ by $H(x,y,t)$, then to solve the equation, it suffices to show the following integral converges,
\begin{eqnarray*}
u(x,t) &=& \int_M H(x,y,t)u_0(y) dy.
\end{eqnarray*}
\begin{equation*}
\abs{\int_{B_t(x,1)}H(x,y,t)u_0(y)dy} \leq Ce^{Ar(x,x_0)},
\end{equation*}
because the integral of $H$ on $B_t(x,1)$ is less than 1 and $u_0$ grows at most exponentially by Theorem \ref{thm_Green}.

\begin{eqnarray*}
\abs{\int_{M\setminus B_t(x,1)} H(x,y,t)u_0(y)dy} &\leq& \int_{M\setminus B_t(x,1)} \frac{C}{V_x(\sqrt{t})}e^{-\frac{r^2(x,y)}{Dt}} \abs{u_0(y)} dy.
\end{eqnarray*}
By volume comparison,
\begin{equation*}
V_x(1)\geq C_1 e^{-A_1r(x,x_0)} V_{x_0}(1)
\end{equation*}
and
\begin{equation*}
V_x(\sqrt{t}) \geq C_2 e^{-A_1r (x,x_0)} \min(1,t^{n/2}).
\end{equation*}
Therefore
\begin{eqnarray*}
\abs{\int_{M\setminus B_t(x,1)} H(x,y,t)u_0(y)dy} &\leq& \int_{M\setminus B_t(x,1)} Ce^{A_1 r(x,x_0)} e^{-\frac{r^2(x,y)}{2DT}}\abs{u_0(y)}dy\\
&\leq&\int_{M\setminus B_t(x,1)} Ce^{A_2 r(x,x_0)}e^{A r(x,y)} e^{-\frac{r^2(x,y)}{2DT}}dy \\
&\leq& Ce^{A_2 r(x,x_0)}.
\end{eqnarray*}

In summary,
\begin{equation*}
\abs{u(x,t)}\leq Ce^{Ar(x,x_0)},
\end{equation*}
where $A$ means a different constant. Volume comparison then implies
\begin{equation*}
	\int_0^T \int_M \exp(-a r^2(x,x_0)) u^2(x,t) dV_t<\infty.
\end{equation*}

For estimates on derivatives, note first that $e^t u(x,t)$ is a solution of heat equation (with evolving metric)
\begin{equation*}
	\pfrac{u}{t}=\triangle u
\end{equation*}
with initial value $u_0$. Since we allow constants depend on $T$, it's equivalent to prove estimates for $e^t u(x,t)$. Therefore, from now on, to the end of this proof, we assume $u(x,t)$ is a solution of heat equation. Then
\begin{equation}
	(\triangle-\pfrac{}{t})u^2=2\abs{\nabla u}.
	\label{eqn_11}
\end{equation}

	Assume that $\varphi: \Real^+\to \Real^+$ satisfies

	1) $\varphi(x)=1$ for $x\leq 1$;

	2) $\varphi(x)=0$ for $x\geq 2$.

	Choose the cut-off function $\varphi(\frac{r(x,x_0)}{R})(R>1)$. Multiplying this to the equation (\ref{eqn_11}) and integrate,
	\begin{equation}
		\int_0^T \int_M\varphi(\frac{r(x,x_0)}{R})\abs{\nabla u}^2 dV_t dt \leq \int_0^T \int_M \varphi(\frac{r(x,x_0)}{R})(\triangle -\pfrac{}{t})u^2 dV_t dt.
		\label{eqn_22}
	\end{equation}
	\begin{eqnarray*}
		\triangle \varphi(\frac{r}{R})&=& \mbox{div} (\varphi^\prime(\frac{r}{R})\frac{1}{R}\nabla r) \\
		&=& \varphi^{\prime\prime}(\frac{r}{R})\frac{\abs{\nabla r}^2}{R^2} +\varphi^\prime(\frac{r}{R})\frac{1}{R}\triangle r.
	\end{eqnarray*}
	By definition of $\varphi$, we know $\varphi(\frac{r}{R})$ vanishes unless $R\leq r(x,x_0)\leq 2R$. Laplacian comparison implies (curvature is bounded from below $-k$)
	\begin{equation*}
		\triangle r\leq (n-1)\sqrt{k}\mbox{coth}(\sqrt{k}r)\leq C.
	\end{equation*}

	Therefore, 
	\begin{equation}
		\label{eqn_33}
		\int_0^T \int_M \varphi \triangle u^2 dV_t \leq C\int_0^T \int_{B(x_0,2R)}u^2 dV_t.
	\end{equation}

	Let $dV_t=e^F dV_0$,
	\begin{eqnarray*}
		\int_0^T \int_M \pfrac{}{t}(\varphi u^2) e^F dV_0 dt &\geq& \int_M \int_0^T \pfrac{}{t}(\varphi u^2 e^F) dt dV_0 - \int_0^T \int_M C \varphi u^2 dV_t dt \\
		&=& \int_M \varphi u^2(x,T) dV_T -\int_M \varphi u^2(x,0)dV_0 -C \int_0^T \int_M \varphi u^2 dV_t dt \\
		&\geq& -\int_M \varphi u_0^2(x) dV_0 -C \int_0^T \int_M \varphi u^2 dV_t dt.
	\end{eqnarray*}
	Here we have used the fact that $\pfrac{e^F}{t}$ is bounded.
	Combined with equation (\ref{eqn_22}) and (\ref{eqn_33}),
	\begin{equation*}
		\int_0^T \int_{B(x_0,R)} \abs{\nabla u}^2 dV_t dt \leq C\int_0^T \int_{B(x_0,2R)} u^2 dV_t dt + \int_{B(x_0,2R)} u_0^2(x) dV_0.
	\end{equation*}
	From here it's easy to see the type of estimate in Theorem \ref{thm_maximum}. For $\triangle u$, it suffices to consider $\abs{\nabla^2 u}$. The Bochner formula in this case is (remember we have assumed that $u$ is a solution of the heat equation),
	\begin{equation*}
		(\triangle -\pfrac{}{t})\abs{\nabla u}^2 = 2\abs{\nabla^2 u}^2 - \abs{\nabla u}^2.
	\end{equation*}
	The same argument as before works for $\abs{\nabla^2 u}$.
\end{proof}

\begin{lem}For $t\in [0,T_{max})$,
\begin{equation*} 
	\triangle_t u(x,t)=R(x,t)+1.
\end{equation*}
\end{lem}

\begin{proof}
	We know for $t=0$ it's true.
	Calculation shows 
	\begin{eqnarray*}
		\pfrac{}{t}(\triangle_t u-R(t)-1) &=& (R+1)\triangle_t u+\triangle_t(\triangle_t u-u) -\triangle_t R -R(R+1) \\
		&=& \triangle_t (\triangle_t u-R(t)-1) +R(\triangle_t u -R -1)
	\end{eqnarray*}
	By previous lemma, we have growth estimate for $\triangle_t u-R(t)-1$. If $\triangle_t u-R-1\geq 0$, then
	\begin{equation*}
		(\pfrac{}{t}-\triangle_t)(\triangle_t u-R(t)-1)\leq C(\triangle_t u-R-1).
	\end{equation*}
	If $\triangle_t u-R-1\leq 0$, then
	\begin{equation*}
		(\pfrac{}{t}-\triangle_t)(\triangle_t u-R(t)-1)\geq C(\triangle_t u-R-1).
	\end{equation*}
	Apply maximum principle for $\triangle_t u -R -1$, which is zero at $t=0$. We know it's zero forever. 
\end{proof}

\section{Proof of the main theorem and the corollary}

Assume we have a surface satisfying the assumptions of Theorem \ref{thm_main}. Short time existence is known, see \cite{Shi}. The long time existence and convergence follows exactly by an argument of Hamilton in \cite{Ha}. For completeness, we outline the steps.

Solve Poisson equations $\triangle_t u(x,t)=R(x,t)+1$ as we did. Consider the evolution equation for $H=R+1+\abs{\nabla u}^2$,
\begin{equation*}
\pfrac{}{t}H=\triangle H-2\abs{M}^2 -H,
\end{equation*}
where $M=\nabla\nabla u - \frac{1}{2}\triangle f\cdot g$.
Since we have growth estimate for $H$, maximum principle says
\begin{equation*}
R+1\leq H\leq Ce^{-t}.
\end{equation*}
Therefore, after some time $R$ will be negative everywhere.
Applying maximum principle again to the evolution equation of scalar curvature
\begin{equation*}
\pfrac{}{t}R=\triangle R+R(R+1)
\end{equation*}
will prove Theorem \ref{thm_main}.

Next, we discuss the application of the above theorem to Uniformization theorem. Let $S$ be a compact Riemann surface. Let $p_1,\cdots,p_k$ be $k$ different points in $S$ and $D_1,\cdots,D_l$ be $l$ domains on $S$ such that all of them are disjoint and $D_i$ is diffeomorphic to disk. Denote $S\setminus \cup_i D_i \setminus \set{p_1,\cdots,p_k}$ by $M$. The aim is to show there exists a complete hyperbolic metric on $M$ compatible with the conformal structure.

The approach is to construct an initial metric $g_0$ on $M$ compatible with the conformal structure so that the normalized Ricci flow starting from $g_0$ will converge to a hyperbolic metric. Assume there is metric $h$ in the given conformal class of $S$. Note that $h$ is incomplete as a metric on $M$.

For $p_i$, there is an isothermal coordinate $(x,y)$ around $p_i$. By a conformal change of $h$, one can ask $g_0$ to be
\begin{equation*}
	g_0=\frac{8}{(x^2+y^2)\log^2 (x^2+y^2)}(dx^2+dy^2)
\end{equation*}
in a small neighborhood $U_i$ of $p_i$.
\begin{rem}
	This is called hyperbolic cusp metric in \cite{Se} and it has scalar curvature $-1$.
\end{rem}

For $D_j$, let $r$ be the distence to $\partial D_j$ on $M$ with respect to $h$. Let $V_j$ be a neighborhood of $\partial D_j$ in $M$. Let $(r,\theta)$ be the Fermi coordinates for  $\partial D_j$ so that
\begin{equation*}
	h_0=dr^2+A(r,\theta)d\theta^2.
\end{equation*}
We will find $\rho=\rho(r,\theta)$ such that

1) $\rho=0$ on $\partial D_j$;

2) $d\rho\ne 0$ on $\partial D_j$;

3)
\begin{equation*}
	g_0=\frac{1}{\rho^2}h
\end{equation*}
is asymptoticly hyperbolic in high order. Let $K$ and $K_0$ be the Gaussian curvature of $h$ and $g_0$ respectively. We have the formula,
\begin{equation*}
	K_0=\rho^2(\triangle_h \log \rho+K).	
\end{equation*}
In order that $K_0=-1$,
\begin{equation*}
	1-\abs{\nabla \rho}^2 +\rho \triangle\rho+\rho^2 K=0.
\end{equation*}
In terms of $r$ and $\theta$,
\begin{equation*}
	\abs{\nabla \rho}^2=(\pfrac{\rho}{r})^2 +A^{-1}(r,\theta) (\pfrac{\rho}{\theta})^2
\end{equation*}
and
\begin{equation*}
	\triangle \rho = \frac{\partial^2 \rho}{\partial r^2} +B\frac{\partial \rho}{\partial r} + C\pfrac{\rho}{\theta} +D \frac{\partial^2 \rho}{\partial \theta^2}.
\end{equation*}
Here $A$, $B$, $C$ and $D$ are smooth functions of $r$ and $\theta$. The equation now becomes
\begin{equation}
	\rho \frac{\partial^2 \rho}{\partial r^2}+B\rho \pfrac{\rho}{r} +C\rho \pfrac{\rho}{\theta} +D\rho \frac{\partial^2 \rho}{\partial \theta^2} +1-(\pfrac{\rho}{r})^2 -A^{-1} (\pfrac{\rho}{\theta})^2 +\rho^2 K=0.
	\label{eqn_aa}
\end{equation}
If equation (\ref{eqn_aa}) is true at $r=0$, then
\begin{equation}
	\label{eqn_bb}
	\pfrac{\rho}{r}(r,\theta)=1.
\end{equation}
Here we used that fact that $\rho>0$.

Set $\eta(r,\theta)=\frac{\rho}{r}$. Equation (\ref{eqn_bb}) implies $\eta(0,\theta)=1$. Equation (\ref{eqn_aa}) becomes
\begin{equation*}
	\begin{array}{l}
		r^2\frac{\partial^2 \eta}{\partial r^2} + Br\eta + Br^2 \pfrac{\eta}{r} +Cr^2 \pfrac{\eta}{\theta} +D r^2\frac{\partial^2 \eta}{\partial \theta^2} \\
		+\frac{1-\eta^2}{\eta} -\frac{r^2}{\eta}(\pfrac{\eta}{r})^2 -A^{-1} \frac{r^2}{\eta}(\pfrac{\eta}{\theta})^2 +\eta r^2 K=0
	\end{array}
\end{equation*}
For the convinience of formal calculation, this equation is rewritten as
\begin{equation}
	(D^2-D-2)\eta+F[r,\eta]=0,
	\label{eqn_cc}
\end{equation}
where $D=r\pfrac{}{r}$ and
\begin{eqnarray*}
	F[r,\eta]&=& Br\eta+Br^2 \pfrac{\eta}{r}+Cr^2 \pfrac{\eta}{\theta} + Dr^2 \frac{\partial^2 \eta}{\partial \theta^2} +\frac{(1-\eta)^2}{\eta}+2 \\
	&&-\frac{1}{\eta}(r\pfrac{\eta}{r})^2 -A^{-1}\frac{r^2}{\eta}(\pfrac{\eta}{\theta})^2 +\eta r^2 K.
\end{eqnarray*}

Equation (\ref{eqn_cc}) is a very typical form of Fuchsian type PDE. Formal solutions of this kind of equation has been discussed many times. For example, Kichenassamy \cite{Kich} and Yin \cite{Yin}. We will only outline the main steps here, for details see \cite{Kich} and \cite{Yin}.

Consider formal solution with the following expansion,
\begin{equation}
	\eta(r,\theta)=1+\sum_{i=1}^\infty\sum_{j=0}^{i} a_{ij}(\theta) r^i(\log r)^j.
	\label{eqn_expansion}
\end{equation}
We will call the sum $\sum_{j=0}^{i}a_{ij}r^i(\log r)^j$ the $i$-level of the expansion. Note that $D$ maps $i$-level to $i$-level. Details on formal calculation could be find in \cite{Kich} and \cite{Yin}. A common feature of all terms in $F[r,\eta]$, which is crutial in obtaining a formal solution, is that the $k$-level of $F[r,\eta]$ could be calculated with knowledge of only $l$-level of $\eta$ with $l<k$. For example, consider $(1-\eta)^2/\eta$. It's the multiplication of three formal series, two $1-\eta$ and $1/\eta$. In order the $k$-level of $\eta$ appears in the $k$-level of $(1-\eta)^2/\eta$, the only possibility is that two of the three series contribute zero level and one $k$-level. However, the zero level of $1-\eta$ vanishes.

The only thing we need is that there exists a formal solution and furthermore due to Borel's Lemma as in \cite{Yin}, there is an approximate solution so that
\begin{equation*}
	(D^2-D-2) \eta +F[r,\eta]=o(r^k)
\end{equation*}
for any $k$. In terms of $\rho$,
\begin{equation}
	K_0+1=1+\rho^2(\triangle_h \log \rho+K)=o(\rho^k)
	\label{eqn_dd}
\end{equation}
for any $k$. This metric $g_0$ near $\partial D_j$ has Gaussian curvature -1 asymptotically. By a scaling, we assume it has scalar curvature -1 asymptotically.

We construct $g_0$ by doing the above to every point $P_i$ and disk $D_j$. If there is at least one disk removed, we know $M$ is nonparabolic. 
\begin{equation*}
	\int_M \abs{R+1}dV
\end{equation*}
is finite because of equation (\ref{eqn_dd}). Therefore, Theorem \ref{thm_main} proves the Uniformization in this case.

If there is no disk removed, i.e. $M=S\setminus \set{p_1,\cdots,p_k}$ and $M$ has negative Euler number, then it's proved in \cite{Se} that there exists a hyperbolic metric in the conformal class. A large part of \cite{Se} is devoted to solve 
\begin{equation}
	\label{eqn_xx}
	\triangle_{g_0} u= R_{g_0}+1
\end{equation}
with $\abs{\nabla u}<\infty$.

Observe that the above equation is equivalent to
\begin{equation}
	\triangle_h u=\frac{g_0}{h}(R_{g_0}+1).
	\label{eqn_final}
\end{equation}
Since every end of $(M,g_0)$ is a hyperbolic cusp, Gauss-Bonnet theorem says
\begin{equation}
	\int_M R_{g_0}dV_0=2\pi \chi(M)<0.
	\label{eqn_zz}
\end{equation}
There exists a function $f$ of compact support on $M$ such that the volume of $(M,e^fg_0)$ is $-2\pi \chi(M)$, because $(M,g_0)$ has finite volume. Denote $e^fg_0$ by $g_0$, since the infinity is not changed, equation (\ref{eqn_zz}) is still true. Now, the volume of $(M,g_0)$ is $-2\pi\chi(M)$. This implies
\begin{equation*}
	\int_M (R_{g_0}+1)dV_0=0.
\end{equation*}
Therefore
\begin{equation*}
	\int_S \frac{g_0}{h}(R_{g_0}+1)dV_h=0.
\end{equation*}
By construction of $g_0$, we know $\frac{g_0}{h}(R_{g_0}+1)$ is zero near $P_i$. So $\frac{g_0}{h}(R_{g_0}+1)$ is a smooth function on $S$. Therefore, equation (\ref{eqn_final}) is solvable. Since $u$ is a smooth function on compact surface $S$, $u$ has bounded gradient with respect to $h$. The relation of $h$ and $g_0$ near $P_i$ is explicit. It's straight forward to check $u$ has bounded gradient as a function of $(M,g_0)$. This simplifies the proof in \cite{Se}. 

\begin{rem}

	In the case that there is at least one disk removed, by construction of $g_0$, $R_{g_0}+1$ vanishes at high order near $\partial D_j$. Then one can extend the definition $\frac{g_0}{h}(R_{g_0}+1)$ to $S$ so that
\begin{equation*}
	\int_S \frac{g_0}{h}(R_{g_0}+1)dV_h=0.
\end{equation*}
The rest is the same as in the previous case.

	This method of solving Poisson equation depends on the conformal structure of $M$, therefore Theorem \ref{thm_Ni} and Theorem \ref{thm_main} are not coverd by the above discussion.
\end{rem}

\end{document}